\documentclass[12pt]{amsart}
\usepackage{amsfonts,amsmath,amscd}
\usepackage{epsfig,graphics}
\usepackage{amssymb}
\usepackage{amscd}
\usepackage{hyperref}

\newtheorem{lemma}{Lemma}[section]

\newcommand{\BN}{{\mathbb{N}}}
\newcommand{\BR}{{\mathbb{R}}}
\newcommand{\BC}{{\mathbb{C}}}

\newcommand{\ga}{\alpha}

\newcommand{\ff}{{\frak{f}}}

\def\De{\Delta}

\newcommand{\vol}{\text{vol}}
\newcommand{\diam}{\text{diam}}

\newcommand{\spn}{\text{span}}
\newcommand{\ad}{\text{ad}}

\newcommand{\End}{\text{End}}

\newcommand{\Hom}{\text{Hom}}
\newcommand{\VN}{\text{VN}}

\newcommand{\Meas}{\text{Meas}}

\def\N{{\cal N}}

\newcommand{\half}{\frac{1}{2}}

\newcommand{\form}[2]{\langle #1,#2 \rangle}

\newtheorem{prop}{Proposition}[section]
\newtheorem{thm}{Theorem}

\newtheorem{cor}[prop]{Corollary}
\newtheorem*{conj}{Conjecture}
\theoremstyle{definition}

\newtheorem{defn}[prop]{Definition}

\newtheorem*{prob}{Problem}

\newcommand{\metric}{\sigma}

\begin{document}

\title[Boundary representations]{Boundary unitary representations - irreducibility and rigidity}
\author{U. Bader \& R. Muchnik}
\thanks{U.B was partially supported by the ISF grant 704/08 and GIF grant 2191-1826.6/2007.}
\thanks{R.M was partially supported by the NSF Postdoctoral Research Fellowship DMS-0202457.}
\maketitle

\begin{abstract}
Let $M$ be compact negatively curved manifold, $\Gamma =\pi_1(M)$ and $\tilde{M}$ be its universal cover.  Denote by $B =\partial \tilde{M}$ the geodesic boundary of $\tilde{M}$ and by $\nu$ the Patterson-Sullivan measure on $X$.  In this note we prove that the associated unitary representation of $\Gamma$ on $L^2(B,\nu)$ is irreducible.  We also establish a new rigidity phenomenon: we show that some of the geometry of $M$, namely its marked length spectrum, is reflected in this $L^2$-representations.

\end{abstract}

\section{Introduction}

Quite often ergodic properties of a measurable dynamical system are reflected in its associated quasi-regular unitary representation.
This is the case, for example, with the notions of {\em ergodicity, mixing} and {\em weak-mixing}
associated to (probability) measure preserving group actions.
They all have an equivalent representation theoretic (or spectral) formulation.
The property of irreducibility of the quasi-regular representation, however, can not occur for probability measure preserving actions,
as the constants always form an invariant line.

\begin{prob}
Given a discrete group can one "classify" its measure class preserving (measurable) actions for which the quasi-regular
representations are irreducible?
\end{prob}

In \cite{Mackey}[\S3.5, Corollary 7] Mackey gives a full solution to the problem when one restricts himself to
(infinite) {\em homogeneous} actions
(see also \cite{Burger-Harpe}).
The general case is still not well understood.
The only existing examples (to the best of the author's knowledge) for measurable non-homogeneous group actions
with irreducible quasi-invariant representations were given in two cases:

\begin{itemize}
\item
Actions of free groups on their boundaries (see \cite{FT-P-book},\cite{FT-Steger}).
\item
Actions of lattices of Lie-groups (or algebraic-groups) on their Furstenberg boundaries (see \cite{Cowling-Steger}, \cite{Bekka-Cowling}).
\end{itemize}

We emphasize that in the above mentioned results one considers boundary
actions of the group in question, and asks whether there is a
general phenomenon yet to be discovered being hinted at by this coincidence.
In particular, we state the following conjecture:

\begin{conj} \label{conj}
For a locally compact group $G$ and a spread-out probability measure $\mu$ on $G$,
the quasi-regular representation associated to a $\mu$-boundary of $G$ is irreducible.
\end{conj}

In this paper we prove this conjecture for the action of the fundamental group of a compact negatively curved manifold on its boundary endowed with the Paterson-Sullivan measure class.
\footnote{Note that free groups are the fundamental groups of space of loops and uniform lattices in rank one Lie groups are fundamental groups of the quotient space.}
We call the associated quasi-regular representation a {\em boundary representation}.
In particular, we prove:

\begin{thm}[Irreducibility] \label{irr}
Every boundary representation is irreducible.
\end{thm}

For every negatively curved manifold,
there is a naturally attached class function of its fundamental group,
called the {\em marked length spectrum}.
The Irreducibility Theorem gives rise to a rigidity statement based on the marked length spectrum.
The following theorem suggests that we view this class function as a
{\em character} associated to the boundary representation.

\begin{thm}[Rigidity] \label{rigidity}
Two boundary representations are unitary equivalent if and only if the
two associated negatively curved manifolds give rise to proportional marked length spectra.
\end{thm}

\subsection{Structure of the paper:}
In the next section we introduce the notations that are used in the paper, and recall some standard results and constructions related to negatively curved manifolds.
We will also state the important Theorems~\ref{main} and \ref{operators} that are needed for the proof of Theorem~\ref{irr}.

Section~3 is dedicated to proofs of our preparatory lemmas,
in particular, there we introduce the "chopped" functions
($\overline{(q|\cdot)}$ and its derivatives), which definitions are motivated by
\cite{Connell-Muchnik-GAFA}.
Section~4 is devoted to the proof of the uniform boundedness of a certain family of functions,
which later on (section~6) will be translated into the pre-compactness of a certain
family of measures.
In section~5 we use a geodesic counting result of Margulis to compute the asymptotic decay
of some matrix coefficients.
In section~6 we combine the results of sections 4 and 5, into proofs
of the theorems stated in section~2, and then show how the latter imply
Theorem~\ref{irr} and Theorem~\ref{rigidity}.

We add three appendices to this short paper.
The first two - concerning with metric measure spaces and functional analysis -
are devoted to the recollection of some basic definitions and proofs of
easy results of general nature, to be used in the body of the paper,
without breaking the principle development of the paper.
The third appendix concerns with a certain
(not very well known)
counting geodesics theorem of Margulis.

\medskip
{\bf Acknowledgments.}
We wish to express our gratitude to Alex Eskin who should have been a third author for the manuscript.  We also would like to thank Fran\c{c}ois Ledrappier who contributed tremendously by sharing his observation about the geometric meaning of the equidistribution result we needed, and pointing out Margulis' thesis as a reference.  Also we wish to express our thanks to Vadim Kaimanovich for pointing us toward the work of Figa-Talamanca and Alex Furman for explaining the various aspects of the marked length spectrum and rigidity results.

We express our appreciation to the following mathematicians for their help:
Chris Connell - for many helpful discussions about negative curvature,
Benson Farb - for his support and the reference he supplied,
Alex Gorodnik - for explaining his equidistribution results and
Nicolas Monod - for various discussions of operator algebras.
We thank the anonymous referee for a careful reading and some helpful suggestions.

\newcommand{\T}{\mathbb{T}}
\newcommand{\E}{\mathbb{E}}
\newcommand{\C}{\mathbb{C}}
\newcommand{\R}{\mathbb{R}}
\renewcommand{\H}{\mathbb{H}}
\newcommand{\HH}{\mathbb{H}}
\newcommand{\Z}{\mathbb{Z}}
\newcommand{\PP}{\mathbb{P}}
\newcommand{\Zp}{\mathbb{Z}_p}
\newcommand{\Q}{\mathbb{Q}}

\newcommand{\eps}{\epsilon }

\newcommand{\set}[1]{\left\{ #1 \right\}}
\newcommand{\op}[1]{\operatorname{#1}}
\newcommand{\Ad}{\operatorname{Ad}}
\newcommand{\Id}{\operatorname{Id}}
\newcommand{\tr}{\operatorname{trace}}
\newcommand{\trace}{\operatorname{trace}}
\newcommand{\supp}{\operatorname{supp}}
\newcommand{\Lip}{\operatorname{Lip}}
\newcommand{\myspan}{\operatorname{span}}
\newcommand{\hrank}{\operatorname{h-rank}}
\newcommand{\htop}{\operatorname{h}}
\newcommand{\Vol}{\operatorname{Vol}}
\newcommand{\Jac}{\operatorname{Jac}}
\newcommand{\Isom}{\operatorname{Isom}}

\newcommand{\A}{\mathsf A }
\newcommand{\G}{\mathsf G }
\newcommand{\Sgp}{\mathsf S }

\newcommand{\meas}{\mathcal{M} }
\newcommand{\D}{\partial }
\newcommand{\M}{\tilde{M}}
\newcommand{\X}{\tilde{X}}
\newcommand{\SM}{SM}
\newcommand{\dM}{\partial\M }
\newcommand{\dN}{\partial\N }

\newcommand{\union}{\cup}
\newcommand{\intersect}{\cap}
\newcommand{\semidirect}{\ltimes}
\newcommand{\directsum}{\oplus}
\newcommand{\grad}{\nabla}
\newcommand{\til}{\widetilde}
\renewcommand{\hat}{\widehat}
\newcommand{\rest}[2]{ #1_{\downharpoonleft_{#2}} }
\newcommand{\dvol}{dg}
\newcommand{\del}{\partial }
\newcommand{\degree}{^\circ}
\newcommand{\degrees}{^\circ}
\newcommand{\goto}[1]{\stackrel{#1}{\longrightarrow}}
\newcommand{\ident}{\equiv }
\newcommand{\comp}{\circ }
\newcommand{\of}{\circ }
\providecommand{\to}{\rightarrow }
\newcommand{\tensor}{\otimes }
\renewcommand{\bar}{\overline}
\newcommand{\minus}{\tiny{-}}
\newcommand{\abs}[1]{\left\lvert #1 \right\rvert }
\newcommand{\inner}[1]{\left\langle #1 \right\rangle }
\newcommand{\dd}[1][t]{\frac{d}{d #1} }
\newcommand{\dds}[1][t]{\frac{d^2}{d #1}^2 }
\newcommand{\omicron}{\circ}
\newcommand{\m}{\frac{d\gamma_* \mu}{d\mu}}
\newcommand{\nin}{\noindent}
\newcommand{\conv}{\star}
\newcommand{\group}{\Gamma}

\renewcommand{\(}{\left(}
\renewcommand{\)}{\right)}
\newcommand{\numeq}[1]{\begin{align}\begin{split} #1
\end{split}\end{align}}

\def\[#1\]{\begin{align*}\begin{split} #1 \end{split}\end{align*} }
\def\[[#1\]]{\begin{align}\begin{split} #1 \end{split}\end{align} }

\def\ssu{\subset}
\def\sss{\supset}
\def\sm{\setminus}
\def\pr{\prime}
\def\zz{ \mathbb Z }
\def\cc{ \mathbb C}
\def\nn{ \mathbb N}
\def\rr{\mathbb R}
\def\ff{ \mathbb F}
\def\qq{ \mathbb Q}
\def\T{ \mathbb T}

\def\De{ \Delta}
\def\Ga{ \Gamma}
\def\ga{\gamma}
\def\om{ \omega}
\def\Om{ \Omega}
\def\Si{ \Sigma}
\def\si{ \sigma}
\def\la{ \lambda}
\def\La{ \Lambda}
\def\til{ \tilde}
\def\<{\langle}
\def\>{\rangle}
\def\pr{\prime}
\def\pa{\partial}
\def\ep{\epsilon}

\newcommand{\marginlabel}[1]
  {\mbox{}\marginpar{\raggedleft\hspace{0pt}#1}}


\section{Preliminaries}

In this section we set the notation, and state some facts, to be used throughout the paper.
We will resist the temptation of stating things in a greater generality than needed.
A good reference for Patterson-Sullivan theory for CAT(-1) space is \cite{BM-CAT(-1)}.

Let $(M,g)$ be a (strictly) negatively curved compact Riemannian manifold.
Up to a rescaling of the metric we can and will assume that $M$ is a CAT(-1) metric space.
Denote the diameter of $M$ by $R$.
We set $\Gamma=\pi_1(M)$ to be the fundamental group of $M$.
Denote by $X$ the universal cover of $M$, 
and by $pr:X \rightarrow M$ the obvious projection map.
The space $X$ is endowed with the lifted Riemannian metric $\til{g}$.
We denote the associated metric of $X$ by $d$.
The ball in $X$ of radius $t\geq 0$, centered at $p\in X$ is denoted $X(p,t)$.
Thus, for example, we have $\Gamma\cdot X(p,R)=X$.
It is known that $\Gamma$ acts isometrically, properly discontinuously and co-compactly on $X$; with the quotient being $\Gamma \backslash X= M$.
Every point $p\in X$ gives rise to length function on $\Gamma$,
defined by $|\gamma|_p=d(p,\gamma p)$.
The {\em Gromov product} on $X$ is defined for every triple$p,q,r \in X$ by
$$(p|q)_r= \frac{1}{2}\left[d(p,r)+d(q,r)-d(p,q)\right],$$
It is well known that $(X,d)$ is a {\em $\delta$-hyperbolic space}.
That is, there exists some $\delta\geq 0$, such that for all $p,q,r,w\in X$.
\begin{equation}
\tag{hyp}
(p| q)_r \geq \min((p| w)_r, (q| w)_r ) -\delta
\end{equation}

Consider the quotient space of complex valued continuous functions on $X$
modulo the constant functions, $C(X)/\C$.
This space is endowed with the quotient Frechet structure, coming from the topology of
uniform convergence on compact sets in $C(X)$.
The map
$$ X\rightarrow C(X)/\C, \quad p\mapsto [d(p,\cdot)] $$
is a homeomorphism on its image, which is relatively compact.
The closure $\bar{X}$ is called {\em the horofunctions compactification} of $X$.
We set $B=\pa X$.
This is a topological sphere.
The $\Gamma$ action on $X$ extends continuously to $\bar{X}$,
with $B$ being the unique minimal closed invariant subset.
The Gromov-product $(\cdot|\cdot)_p$ continuously extends to $\bar{X}\times\bar{X} $
(for a fixed $p\in X$),
and the inequality (hyp) stays valid.
For a fixed a point $p\in X$, the function
$$\metric_p:B\times B\rightarrow \R, \quad \metric_{p}(b,c) = e^{-(b | c)_p}$$
is a metric,
denoted {\em the Busemann metric}, on $B$ (see \cite[Chapter 2]{Bourdon}).
The ball in $B$ with respect to $\metric_p$, of radius $t\geq 0$, centered at $b\in B$ will be denoted
$B_p(b,t)$.
For every two distinct points $r\in X$, and $q\in \bar{X}$ denote by
$\ell_{r,q}:[-\infty,\infty]\rightarrow \bar{X}$ the unique unit length geodesic passing through $r$ and $q$,
and satisfying $\ell_{r,q}(0)=p$ and $\ell_{r,q}(d(r,q)) =q$.
Denote $\ell_{r,q}(\infty)$ by $z_{r}^q$ (observe that if $q \in B$ then $z_r^q=q$).
We denote $B_p(q)=B_p(z_p^q,e^{-d(p,q)})$,
and (for completeness) $B_p(p)=B$.
When $p$ is fixed and $q$ varies the map $q\mapsto B_p(q)$ is a surjective continuous map
from $X$ to the collection of balls of $(B,\metric)$, taken with the Hausdorff distance.
The quantity
$$ \beta_{b}(p,q) = d(p,q)- 2 (q | b)_p,$$
is defined for all $p,q \in X$ and $b \in B$.
It is called {\em the Busemann cocycle},
and satisfies the {\em cocycle equation}
$$\mbox{for every } p,q,r \in X,~b \in B,\quad
\beta_b(p,q)+\beta_b(q,r)= \beta_b (p,r).$$

Denote by $\eta$ the {\em critical exponent} associated to $(M,g)$, that is
(for a fixed, yet insignificant, $p\in X$)
$$ \eta=\lim_{t\rightarrow\infty} \frac{\log\vol(X(p,t))}{t}. $$
The critical exponent is positive and finite.
It is known that $\eta$ equals the {\em Hausdorff dimension} of the metric space $(B,\metric_p)$.
The associated {\em Hausdorff measure}
is denoted $\nu_p$.
The collection of measures $\nu_p$, $p\in X$ is
an $\eta$-conformal density for $\Gamma$,
i.e.,the associated map
$\nu:X\rightarrow Meas(B)$ is $\Gamma$-equivariant,
and for every $p,q\in X$,
$$ \frac{d\nu_q}{d\nu_p}(b)=e^{-\eta\beta_b(p,q)}. $$
where $\Meas(B)$ denotes the space of radon measures on $B$.

It is known that $\nu$ is the unique (up to a constant multiplicative factor)
$\eta$-conformal density.
In particular it coincides (up to a scale) with the {\em Patterson-Sullivan measures}.
Therefore, we choose the representation such that $\nu_p(B)=1$ for some fixed point $p \in X$.

Denote $H=L^2(B,\nu_p)$ the space of complex valued square integrable functions on $B$.
The vector space $H$ is independent of $p$.
Endow $H$ with the inner product
$$\mbox{for every } u,v\in H ,\quad \form{u}{v}_p=\int_B u(b)\overline{v(b)}d\nu_p(b).$$
We denote by $\mathcal{U}_p(H)$ the group of unitary operators of $(H,\form{\cdot}{\cdot}_p)$.
The map
$$\rho_p:\Gamma\rightarrow \mathcal{U}_p(H),\quad
\rho_p(\gamma)v(b)=e^{-\half\eta\beta_b(p,\gamma p)}\cdot v(\gamma^{-1}b) $$
is a unitary representation of $\Gamma$,
called the {\em quasi-regular representation},
or the {\em boundary representation} associated to $(B,\nu_p)$.
For every two points $p,q\in X$, the boundary representations $\rho_p$ and $\rho_q$ are unitary equivalent
by the intertwining map
$v\mapsto e^{-\half\eta\beta_b(q,p)}\cdot v$.

Set $H^+=\{v\in H~|~v\geq 0\}$. Observe that this is  a closed cone in $H$.
We denote by $\End(H)$ the algebra of bounded operators on $H$
and endow it with the weak operator topology.
We denote by $\End^+(H)$ the cone of all {\em positive bounded operators},
namely all operators $T\in \End(H)$ that satisfy $TH^+\subset H^+$.
We set $\VN(\rho_p)$ to be the closure of $\spn(\rho_p(\Gamma))<\End(H)$,
and $\VN^+(\rho_p)$ to be the closed cone generated by $\rho(\Gamma)$,
that is, the
closure of the cone of linear combinations
with {\em positive coefficients} of elements of $\rho_p(\Gamma)$.
Obviously, $\VN^+(\rho_p) \subset \End^+(H)$.
The following is the main theorem of this paper.

\begin{thm}[main theorem] \label{main}
$\VN^+(\rho_p) = \End^+(H)$.
\end{thm}

We now present a family of group algebra elements which play a central role in our analysis.
Fix $p \in X$.
For every $t>0$ set
$$ \Gamma \supset S_t=\{\ga~|~ d(\ga p,p) \in (t-R,t+R) \}, $$
and to every bounded Borel function on on $f:B\to \C$ we associate the group algebra element
$$ T^f_t
=\frac{1}
{ |S_t|}
\sum_{\ga \in S_t}\frac{f(z_p^{\ga p})}{\form{\rho_p(\gamma) 1}{1}_p} \cdot \ga \in \C \Gamma. $$
The representation $\rho_p$ extends linearly to the group algebra $\C\Gamma$.
We set $T_t(f) =  T^f_t$.
Thus, $\rho_p \circ T_t$ might be seen as the $\End(H)$-valued measure on $B$.
The space of $\End(H)$-valued measures carries a natural weak$^*$-topology
(see appendix~\ref{app-op}, and in particular, lemma~\ref{op-valued-measure}).
An example of an $\End(H)$ valued measure is given by the measure $m_p$ defined by
$$ C(B) \ni f \mapsto m_p(f), \quad \mbox{where, for } v \in H, \quad
m_p(f)v = \int_B v d\nu_p \cdot f. $$

\begin{thm}[measure convergence] \label{operators}
With respect to the weak$^*$ topology, we have
$$ \lim_{t\rightarrow\infty} \rho_p \circ T_t= m_p. $$
In particular, for $f\in C(B)$ and $g,h \in L^2 (B,\nu_p)$ we have
$$\lim_{t\rightarrow\infty} \frac{1}
{|S_t|}
\sum_{\ga \in S_t}\frac{f(z_p^{\ga p})}{\form{\rho_p(\gamma) 1}{1}_p}\< \rho_p(\ga)g,h\>_p = \left(\int_B f(b)g(b) d\nu_p(b)\right)\left(\int_B h(b)d\nu_p(b)\right).$$
\end{thm}

\section{Some Lemmas}

In this section we start the analysis of some functions
on the boundary, $B$.
Before starting we will establish our point of view.
{\bf We fix once, and for the rest of the paper, a point $p\in X$}.
We set $|q|=d(p,q)$.

We now proceed by defining some family of important functions on the boundary,
and their "chopped" analogs.
Set $\lambda^q(b)=e^{-\half\eta\beta_b(p,q)}$,
and $\|\lambda^q\|_1=\form{\lambda^q}{1}$.
Sometimes we use $\lambda_{\gamma}=\lambda^{\gamma p}=\rho(\gamma)1$
and $\|\la_{\gamma}\|_1=\|\lambda^{\gamma p}\|_1=\form{\rho(\gamma)1}{1}$.
Define the functions $\overline{(q|\cdot)}$ on $B$ by
$$ \overline{(q|b)} = \min\{(z_p^q|b),|q|\}, $$
and accordingly,
$$ \bar{\beta}_b(p,q)=|q|-2 \overline{(q|b)}, \quad \bar{\lambda}^q=e^{-\half\eta\bar{\beta}_b(p,q)}.$$
The justification for the notation is given by

\begin{lemma} \label{comparison}
For every $p,q\in X$ and for every $b\in B$,
\begin{align*}
& 1) & |\overline{(q|b)}-(q|b)| \leq \delta &\\
& 2) & |\bar{\beta}_b(p,q)-\beta_b(p,q)| \leq 2\delta &\\
& 3) & e^{-\delta \eta}\lambda^q(b) \leq \bar{\lambda}^q(b) \leq e^{\delta \eta }\lambda^q(b) &\\
& 4) & e^{-\delta \eta}\|\lambda^q\|_{1} \leq \|\bar{\la}^q\|_1 \leq e^{\delta \eta}\|\lambda^q\|_{1} &
\end{align*}
\end{lemma}

\begin{proof}
Parts (2),(3) and (4) follow from (1), which we will prove.
By (hyp) and $(z_p^q|q) =|q|$ we have
\begin{equation}
(q|b) \geq \min\{(z_p^q|b), (q|z_p^q)\}-\delta=\min\{(z_p^q|b), |q|\}-\delta=\overline{(q|b)}-\delta.
\end{equation}
On the other hand as $(q|b) \leq |q|$
\begin{equation}
(z_p^q|b) \geq \min\{(z_p^q|q), (q|b)\} -\delta = \min\{|q|,(q|b)\}-\delta = (q|b)-\delta.
\end{equation}
Therefore
\begin{equation}
\overline{(q|b)} = \min\{(z_p^q|b),|q|\} \geq (q|b)-\delta
\end{equation}
This finishes the proof. \end{proof}

For the reader convenience the definition of a {\em regular metric measure space}
and a short discussion of its properties is given in appendix~\ref{app-mms}
(definition~\ref{d:mms}).
We set $\metric=\metric_p$.

\begin{prop}
$(B,\metric,\nu)$ is $\eta$-regular.
\end{prop}

We will denote the associated multiplicative constants by $0<k\leq k'$.
Recall our definition $B_r(q)=B_r(z_r^q,e^{-d(r,q)})$ and set also $B(q)=B_p(q)$.

\begin{cor}[{\ref{int-as-log}}]
For every $q\in X$ we have
$$ k|q| -(k'-k) \leq \int_{B(q)^c} \metric(b,c)^{-\eta} d\nu(c) \leq k' |q| + (k'-k)$$
\end{cor}

Before proving the proposition we need to prove the following lemma.

\begin{lemma}
$\inf\{\nu_q(B(q))~|~q\in X\} >0$.
\end{lemma}

\begin{proof}[Proof of the lemma]
We first claim that for every $\gamma\in \Gamma$,
$$B_{\ga q}(\gamma z_p^q,e^{-\delta}) \subset \gamma B_p(q).$$
Indeed, fix a point $\gamma b \in B_{\ga q}(\gamma z_p^q,e^{-\delta})$,
and observe that
$$ 0=(p|z_p^q)_q \geq \min\{(p|b)_q,(b|z_p^q)_q\} -\delta =
  \min\{(p|b)_q-\delta,(\gamma b|\gamma z_p^q)_{\gamma q} -\delta\}. $$
As $(\gamma b|\gamma z_p^q)_{\gamma q}> \delta$, we get $(p|b)_q \leq \delta$.
Therefore
$$
(b|z_p^q)_p
=|q|+ (b|z_p^q)_q - (p|b)_q
=|q|+ (\gamma b|\gamma z_p^q)_{\gamma q} - (p|b)_q
>|q|+ \delta  - \delta =|q|, $$
which proves the claim.

Recall now that $\Gamma \cdot X(p,R)=X$,
and choose $\gamma\in \Gamma$ such that $\gamma q\in X(p,R)$.
From the definition of the Gromov product we get that for every
$x,y\in \bar{X}$, $(x|y)_{\gamma q} \geq (x|y)_p-R$, hence
$$ B_p(\gamma z_p^q|e^{-(\delta+R)}) \subset B_{\gamma q}(\gamma z_p^q,e^{-\delta}). $$
Since  $\nu$ is $\Gamma$-equivariant, we have
$$ \nu_q(B_p(q))=\nu_{\gamma q}(\gamma B_p(q))
\geq \nu_{\gamma q}(B_{\ga q}(\gamma z_p^q,e^{-\delta})
\geq \nu_{\gamma q}B_p(\gamma z_p^q|e^{-(\delta+R)}). $$
By compactness argument, using lemma~\ref{l:lower-cont} in the appendix
$$ \inf\{\nu_{p'}(B_p(b,e^{-(\delta+R)}))~|~p'\in X(p,R),~b\in B\} >0. $$
This proves the lemma.
\end{proof}

\begin{proof}[Proof of the proposition]
Any ball in $B$ equals $B(q)$ for some $q\in B$.
By the $\eta$-conformality of the measure,
$$ \nu(B(q))=\int_{B(q)} d\nu =\int_{B(q)} e^{\eta\beta_b(p,q)} d\nu_q(b) $$
By lemma~\ref{comparison}(2), on $B(q)$, $|\beta_b(p,q)+|q|| \leq 2\delta$, hence
$$
\int_{B(q)} e^{\eta\beta_b(p,q)} d\nu_q(b) \geq
e^{-\eta(|q|+2\delta)}\cdot\inf\{\nu_q(B(q))~|~q\in X\},
$$
and
$$
\int_{B(q)} e^{\eta\beta_b(p,q)} d\nu_q(b) \leq
e^{-\eta(|q|-2\delta)}\cdot\sup\{\nu_q(B(q))~|~q\in X\}.
$$
Denoting
$$
k=e^{-2 \eta \delta} \cdot \inf\{\nu_q(B(q))~|~q\in X\},
\quad \mbox{and} \quad
k'=e^{2 \eta \delta} \cdot \sup\{\nu_q(B(q))~|~q\in X\}.
$$
We get
$$ k \left(e^{-|q|}\right)^{\eta} \leq \nu(B(q)) \leq k' \left(e^{-|q|}\right)^{\eta} $$
This proves $\eta$-regularity.
\end{proof}

%
%
%

We will also need an estimate of the $L^1$ norm of $\la^q$ and $\bar{\la}^q$.
\begin{prop} \label{p:Lambda-estimation}
There exist constants $0<C<C'$ and $0<T$ such that for every $q\in X(p,T)^c$,
$$ C|q|e^{-\half\eta|q|} \leq \|\bar{\la}^q\|_1 \leq
C' |q|e^{-\half\eta|q|},
$$
and
$$ e^{-\delta}C|q|e^{-\half\eta|q|} \leq \|\lambda^q\|_{1} \leq
e^{\delta}C' |q|e^{-\half\eta|q|}.
$$
\end{prop}

\begin{proof}
First observe that for every $q\in X$, by the definition of $\overline{(q|b)}$,
$$ \int_B e^{\eta \overline{(q|b)}} d\nu(b)=
e^{\eta|q|}\nu(B(q)) +
\int_{B(q)^c} e^{\eta (b|\ell_q)} d\nu(b)
$$
Applying $\eta$-regularity to the first summand of the right hand side,
and corollary~\ref{int-as-log} to the second, we get
$$ k|q| -(k'-2k) \leq \int_B e^{\eta \overline{(q|b)}} d\nu(b) \leq k' |q| + (2k'-k)$$
Thus, picking $C<k$ and $C'>k'$, there exist $T>0$ such that for all $|q|>T$,
$$ C|q| \leq \int_B e^{\eta \overline{(q|b)}} d\nu(b) \leq C' |q|$$
As $\|\bar{\la}^q\|_1=\int_B e^{\eta(\overline{(q|b)}-\half |q|)} d\nu(b)$, we get
$$ C|q|e^{-\half \eta|q|} \leq \|\bar{\lambda}^q\|_1 \leq
C' |q|e^{-\half \eta|q|}.
$$
The second assertion follows by lemma~\ref{comparison}.

\end{proof}

\begin{cor} \label{c:unif-bounded-ratios}
For every $q,q'\in X$ with $|q|,|q'|>T$, the ratio of every two functions among
$\|\lambda^q\|_{1},\|\lambda^{q'}\|_1,\|\bar{\lambda}^q\|_1$ and $\|\bar{\lambda}^{q'}\|_1$ is bounded by a constant dependent on $|q|-|q'|$.
\end{cor}

\section{Uniform boundedness} \label{s:bounded}

This section is devoted to the proof that the functions
$$(\rho\circ T^1_t)(1)=\frac{1}{|S_t|}
\sum_{\ga \in S_t}\frac{1}{\form{\rho(\gamma) 1}{1}} \cdot \rho(\ga)1 =\frac{1}{|S_t|}
\sum_{\ga \in S_t}\frac{\la^{\ga p}}{\|\la^{\ga p}\|_1}$$
are uniformly bounded in $t$ in $L^{\infty}(B,\nu)$.
This is achieved in proposition~\ref{p:bounded}.
Our main tool is lemma~\ref{l:sampling} which, after a minor effort,
we show applies here.

The following lemma computes the radius of the
inside and outside balls of a shape in $X$ which arises as the intersection
of an annulus with a cone based at a boundary ball (see figure).

\begin{lemma} \label{shadow-comparison}
Let $q\in X$ be a point such that $q\neq p$.
Denote  $$q'=\ell_{p,q}(|q|+2\delta+ R),$$
and define the set
$$ Y^q=\{ r\in X~|~z_p^r\in B(q),~ \left||r|-|q'|\right| <  R \}, $$
Then $X(q', R) \subset Y_q^p \subset X(q,4\delta+2R)$.
\end{lemma}
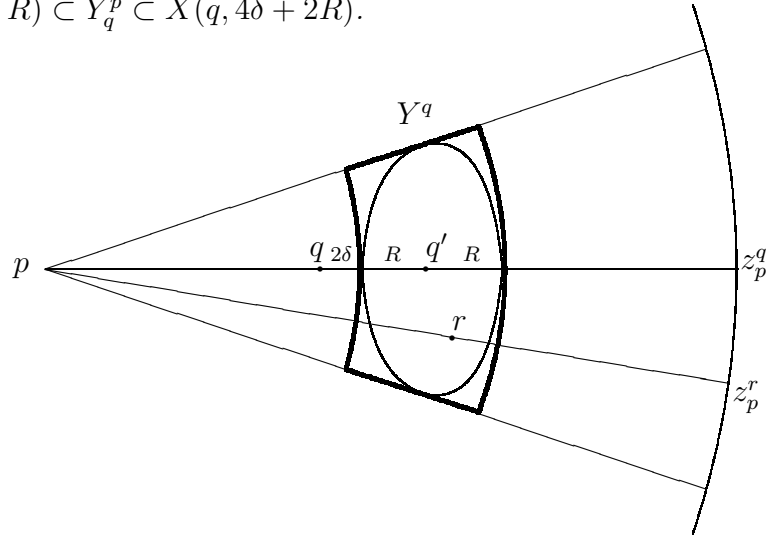
\begin{picture}(300,200)(-350,0)
\put(-300,119){$p$}
\put(-288,120){\line(3,1){250}}
\put(-288,120){\line(3,-1){250}}
\put(-288,120){\line(6,-1){258}}\put(-28,70){$z_p^r$}
\put(-288,120){\line(1,0){262}} \put(-25,120){$z_p^q$}
\put(-144,120){\circle*{2}} \put(-144,125){$q'$}
\put(-184,120){\circle*{2}} \put(-188,125){$q$}
\put(-160,123){\tiny{$R$}}
\put(-130,123){\tiny{$R$}}
\put(-180,123){\tiny{$2\delta$}}
\put(-134,94){\circle*{2}} \put(-134,97){$r$}
\qbezier(-43,220)(-10,120)(-43,20)

\qbezier(-145,73)(-166,80)(-168,120)
\qbezier(-145,73)(-120,67)(-115,120)
\qbezier(-145,167)(-120,173)(-115,120)
\qbezier(-145,167)(-166,160)(-168,120)
\linethickness{.5mm}
\qbezier(-174,157)(-164,120)(-174,83)
\qbezier(-124,174)(-104,120)(-124,66)
\qbezier(-124,174)(-149,166)(-174,158)
\qbezier(-124,66)(-149,74)(-174,82)

\put(-155, 175){$Y^q$}
\end{picture}

\begin{proof}
First, we prove the left hand inclusion.
Fix a point $r\in X(q',R)$.
Obviously, $||r|-|q'|| < R$.
By lemma~\ref{comparison}(1), $$(z_p^q|z_p^r) \geq
\overline{(q'|z_p^r)} \geq (q'|z_p^r)-\delta.$$
As $|r| \geq |q'|-R \geq |q|+2\delta$, by (hyp),
$$ (q'|z_p^r)\geq \min\{(r|z_p^r),(r|q')\}-\delta \geq
\min\{|r|,\half(|r|+|q'|-R)\}-\delta \geq
|q|+\delta.$$
Therefore, $(z_p^q|z_p^r)\geq |q|$ and $z_p^r \in B(q)$.

Now we prove the right hand inclusion.
Fix a point $r\in Y^q$.  Then $z_p^r \in B(q)$.
By Lemma~\ref{comparison}(1),
$(q'|c) \geq |q'|-\delta$.
By (hyp)
$$ (q|r) \geq \min\{(r|z_p^r),(q|z_p^r)\}-\delta
\geq \min\{|r|,|q|-\delta\}-\delta
= |q|-2\delta. $$
Therefore,
$d(q,r)=|q|+|r|-2(q|r) \leq |r|-|q|+4\delta \leq 4\delta +2R$.

%
%


\end{proof}

A consequence of lemma~\ref{shadow-comparison} and the definition of $R$
is that for every $q\in X$,
$$ 1 \leq |\{\gamma\in \Gamma~|~\gamma p \in Y^q \}| \leq m$$
for some fixed integer $m$ (one can take $m=|\{\gamma~|~|\gamma p|<3R+4\delta\}|$).

\begin{defn}
Fix $r>0$.
An $r$-sampling set for $(B,\metric)$, with multiplicity $m\in \BN$,
is a finite set $S$ together with a map $\tilde{\ell}:S\rightarrow B$, $s\mapsto \tilde{\ell}_s$,
such that
$B=\cup_S B(\tilde{\ell}_s, r)$,
and for all $b\in B$, $|\tilde{\ell}^{-1}(B(b, r))| \leq m$.
\end{defn}

\begin{cor} 
For every $t>R+2\delta$,
the set $S_t=\{ \gamma ~|~||\gamma p|-t| \leq R\}$ together with the map $\ell:S_t \rightarrow B$ defined as $\ell(\ga)=z_p^{\ga p}$
is an $e^{-t+R+2\delta}$-sampling set in $B$ with multiplicity $m$.
\end{cor}

\begin{proof}
Fix $b\in B$.  Apply Lemma \ref{shadow-comparison} with  $q=\ell_{p,b}(t-R-2\delta)$ and $q'=\ell_{p,b}(t)$.
\end{proof}

By applying lemma~\ref{l:sampling} we get the following global estimate.

\begin{cor} \label{c:esstimation}
Fix $t>R+2\delta$.
For every $q \in  X$, with $|q|\leq t+R$,
$$ \frac{1}{|S_t|} \sum_{S_t} {\lambda}^q(z_p^{\ga p}) \leq C\|{\lambda}^q\|_1 $$
\end{cor}

\begin{proof} 
Observe that for all $b,c \in B$ such that $(b|c) \geq d(p,q)-2R-2\delta$ (equivalently $\metric(b,c)\leq e^{-d(p,q)+2R+2\delta}$) we
have
\begin{align*}
(b|q) &\geq \min((b|c), (c|q)) -\delta \\ &\geq \min(d(p,q)-2R-2\delta, (c|q)) -\delta = (c|q) -2R-3\delta,
\end{align*}
as $(c|q) \leq d(p,q)$.
Therefore for all $b, c$ such that $\metric(b,c)\leq e^{-t+R+2\delta}\leq e^{-d(p,q)+2R+2\delta}$ we have
\begin{align*}|\log(\la^q(b))-\log(\la^q(c))|&=|\frac{\eta}{2}(\beta_b(p,q)-\beta_c(p,q))|
=|\eta((b|q)-(c|q))|\\ &\leq \eta(2R+3\delta)=L\end{align*}
We use lemma~\ref{l:sampling} to finish the proof. \end{proof}

\begin{prop} \label{p:bounded}
The function (in the variable $t$),
$$ \left\| \frac{1}{|S_t|} \sum_{S_t} \frac{\lambda^{\gamma p}}{\|\lambda^{\gamma p}\|_1} \right\|_{\infty} $$
is bounded on $(R+2\delta,\infty)$.
\end{prop}

\begin{proof} 
Fix $t>R+2\delta$.
Fix $b\in B$ and set $q=\ell_{p,b}(t+R)$.  Thus, $b=z_p^q$.
Observe that for every $y\in X(p,R+t)-X(p,t-R)$,
$$\overline{(y|b)}=\min\{(z_p^y|b),|y|\} \leq
\min\{(z_p^y|b),|q|\} =\min\{(z_p^q|z_p^y), |q|\} =\overline{(q|z_p^y)}.$$
Thus by Lemma \ref{comparison}(1)
$$(y|b) \leq (q|z_p^y)+2\delta.$$
This implies that
$$\beta_b(p,y)=|y|-2(y|b)\geq |q|-2R-2(q|z_p^y)=\beta_{z_p^y}(p,q)-2R.$$
Hence $\la^y(b) \leq \la^q(z_p^y)e^{ \eta R}$.
Since for every $\ga \in S_t$, $\ga p \in X(p,R+t)-X(p,t-R)$ we have
$$ \frac{1}{|S_t|} \sum_{S_t} \frac{\lambda^{\gamma p}(b)}{\|\lambda^{\gamma p}\|_1}
\leq  \frac{e^{\eta R}}{|S_t|} \sum_{S_t} \frac{\lambda^q(z_p^{\ga p})}{\|\lambda^{\gamma p}\|_1}. $$
The proposition now follows by combining corollary~\ref{c:unif-bounded-ratios}
with corollary~\ref{c:esstimation}.
\end{proof}

\section{Analyzing matrix coefficients} \label{s:limits}

In this section the asymptotic of some matrix coefficients of operators
of the form $T_t^{\chi_U}$ (for appropriate subsets $U \subset B$)
is computed, with the aid of appendix~\ref{app-Margulis},
thus providing the main technical tool needed in the proof of theorem~
\ref{operators} given in the next section.
The outcome of this section will be proposition~\ref{p:limits}.
We begin with the more general setting of proposition~\ref{p:limsup}.

Given a positive number $a>0$ and a subset $U \subset B$ we set
$$ U(a)=\{b\in B~|~(b|U) > a\} = \{ b\in B~|~\metric(b,U) < e^{-a}\} $$
Thus, $\cap_{a>0} U(a)=\bar{U}$.

\begin{prop} \label{p:limsup}
Assume we are given a family of elements $\psi_t\in\ell^1(\Gamma)$,
where $t$ is a real positive parameter
such that for every $t$, $\|\psi_t\|_1 \leq 1$, and for every $\ga\in \Gamma$,
$$ \lim_{t\rightarrow \infty} \psi_t(\ga) = 0. $$
Then for every measurable set $V \subset B$, and every $a>0$,
$$ \limsup_{t\rightarrow \infty}
\sum_{\ga\in\Gamma} \psi_t(\ga)\frac{\form{\rho(\ga)1}{\chi_V}}{\form{\rho(\ga)1}{1}}
\leq \limsup_{t\rightarrow \infty}
\sum_{\ga\in \Gamma} \psi_t(\ga)\chi_{V(a)}(z_p^{\ga p}). $$
\end{prop}

For the proof we will state and prove the following lemma.

\begin{lemma}
There exists a constant $C_0$ such that for $q \in X$,
if $\ell_{p,q}(\infty)=z_p^q \notin V(a)$
then
$$\frac{\<\la^q, \chi_V\>}{\|\la^q\|_1} \leq \frac{C_0 e^{a}}{|q|}.$$
In particular if, for some $\ga\in \Gamma$, $\ell_{p,\ga p}(\infty) = z_p^{\ga p} \notin V(a)$
then
$$\frac{\<\rho(\ga)1, \chi_V\>}{\<\rho(\ga)1, 1\>} \leq \frac{C_0 e^{a}}{|\ga p|}.$$
\end{lemma}

\begin{proof} Observe that for $|q|\leq a$ we have
$$\frac{\<\la^q, \chi_V\>}{\|\la^q\|_1} \leq 1\leq \frac{e^a}{a}.$$
For $b \in V(a)$  and $|q|\geq a$ we have
$$\bar{(b|q)} = \min\{(b|z_p^q) , |q|\}\geq  \min\{a , |q|\} =a.$$
Observe that, by lemma \ref{comparison},
\begin{align*}
\form{\la_q}{\chi_V} &\leq \form{e^{\delta \eta} \bar{\la}_q}{\chi_V}
= e^{\delta \eta- \half\eta|q|}\int_V e^{\eta\overline{(q|b)}} d\nu(b)
= e^{\delta \eta-\half\eta|q|}\int_V e^{\eta a} d\nu(b) \\
&\leq e^{\delta+\eta a} \nu(B) e^{-\half\eta |q|}
\end{align*}
We are done by proposition \ref{p:Lambda-estimation}.
\end{proof}

\begin{proof}[Proof of proposition {\ref{p:limsup}}]
Fix a measurable set $V\subset B$, and a positive number $a$.
Fix $t>t_0>0$. Define $\ell:X\to B$ as $\ell(q)=z_p^q$
Decomposing $\Gamma=\Gamma_1 \cup \Gamma_2 \cup \Gamma_3$, using
$\Gamma_1=\{ \ga~|~|\ga p|<t_0\}$, $\Gamma_2=\Gamma_1^c \cap \ell^{-1}(V(a))$
and $\Gamma_3=\Gamma_1^c \cap \ell^{-1}(V(a))^c$,
we get
\begin{align*}
&
\sum_{\Gamma} \psi_t(\ga)\frac{\form{\rho(\ga)1}{\chi_V}}{\form{\rho(\ga)1}{1}} =
\\ &
\sum_{\Gamma_1} \psi_t(\ga)\frac{\form{\rho(\ga)1}{\chi_V}}{\form{\rho(\ga)1}{1}} +
\sum_{\Gamma_2} \psi_t(\ga)\frac{\form{\rho(\ga)1}{\chi_V}}{\form{\rho(\ga)1}{1}} +
\sum_{\Gamma_3} \psi_t(\ga)\frac{\form{\rho(\ga)1}{\chi_V}}{\form{\rho(\ga)1}{1}} \leq
\\ &
\sum_{\Gamma_1} \psi_t(\ga) +
\sum_{\Gamma_2} \psi_t(\ga) +
\sum_{\Gamma_3} \psi_t(\ga) \frac{C_0 e^{a}}{|t_0|}
\end{align*}
Taking limsup of both sides,
the first summand of the right hand side vanishes by assumption, and we get
$$ \limsup_{t\rightarrow \infty}
\sum_{\ga\in\Gamma} \psi_t(\ga)\frac{\form{\rho(\ga)1}{\chi_V}}{\form{\rho(\ga)1}{1}}
\leq \limsup_{t\rightarrow \infty}
\sum_{\ga\in \Gamma} \psi_t(\ga)\chi_{V(a)}(z_p^{\ga p})
+ \frac{C_0 e^{a}}{|t_0|} $$
as $t_0$ was arbitrary, the lemma is proved.
\end{proof}

Specializing to $$\psi_t=T^{\chi_U}_t=\frac{1}
{ |S_t|}
\sum_{\ga \in S_t}\frac{\chi_U(z_p^{\ga p})}{\form{\rho(\gamma) 1}{1}} \cdot \ga$$ we immediately get

\begin{cor}
For measurable subsets $U,V \subset B$, with $\metric(U,V)>0$,
$$ \lim_{t\rightarrow \infty} \form{T_t^{\chi_U}(1)}{\chi_V} =0 $$
\end{cor}

On the other hand, specializing to $$\psi_t=\Check{T}^{\chi_U}_t=\frac{1}
{ |S_t|}
\sum_{\ga \in S_t}\frac{\chi_U(z_p^{\ga p})}{\form{\rho(\gamma) 1}{1}} \cdot \ga^{-1}$$ we deduce

\begin{cor}
For every measurable subsets $U,V\subset B$ and for every $a>0$,
\begin{align*}
&
\limsup_{t\rightarrow \infty} \form{T^{\chi_U}_t(\chi_V)}{1} \leq
\limsup_{t\rightarrow \infty} \frac{1}{|S_t|}
\sum_{\ga\in S_t} \chi_U(z_p^{\ga^{-1} p})
\chi_{V(a)}(z_p^{\ga p})
\end{align*}
\end{cor}

Finally, we get

\begin{prop} \label{p:limits}
For every measurable sets $U,V,W\in B$
such that $\nu(\partial U)=\nu(\partial V)=\nu(\partial W)=0$,
$$
\lim_{t\rightarrow \infty} \form{T^{\chi_U}_t (\chi_V)}{\chi_W} =
\nu(U\cap W)\nu(V) $$
\end{prop}

\begin{proof}
The last two corollaries
readily imply that for every $a,a'>0$,
\begin{align*}
&
\limsup_{t\rightarrow \infty} \form{T^{\chi_U}_t (\chi_V)}{\chi_W}
\leq
\limsup_{t\rightarrow \infty} \frac{1}{|S_t|}
\sum_{\ga\in S_t} \chi_{U\cap W(a)}(z_p^{\ga^{-1} p})
\chi_{V(a')}(z_p^{\ga p})
\end{align*}
By corollary~\ref{c:margulis} in the appendix,
for every $a$ and $a'$ such that
$\nu(\partial W(a))=\nu(\partial V(a'))=0$,
$$ \lim_{t\rightarrow \infty} \frac{1}{|S_t|}
\sum_{\ga\in S_t} \chi_{U\cap W(a)}(z_p^{\ga^{-1}p})
\chi_{V(a')}(z_p^{\ga p}) =
\nu(U\cap W)\nu(V). $$
This condition is valid for all,
but at most countable many values of $a$ and $a'$.
Thus, by taking sequences tending to zero, we get
$$
\limsup_{t\rightarrow \infty} \form{T^{\chi_U} (\chi_V)}{\chi_W} \leq
\nu(U\cap W)\nu(V). $$
From the definition of $T$ it is obvious that $\form{T_t^1(1)}{1}=1$
for every $t>0$.
If we denote $U^1=U,~U^{-1}=B-U$, and similarly for $V$ and $W$,
we get
\begin{align*}
 1 = & \liminf_{t\rightarrow \infty} \form{T^1_t(1)}{1}
\\ \leq &
\liminf_{t\rightarrow \infty} \form{T^{\chi_U}_t (\chi_V)}{\chi_{W}}
+
\sum_{\substack{i,j,k=\pm 1 \\ i,j,k \neq 1,1,1}}
\limsup_{t\rightarrow \infty} \form{T^{\chi_{U^i}}_t (\chi_{V^j})}{\chi_{W^k}}
\\ \leq &
\sum_{i,j,k=\pm 1}
\limsup_{t\rightarrow \infty} \form{T^{\chi_{U^i}}_t (\chi_{V^j})}{\chi_{W^k}}
\leq
\sum_{i,j,k=\pm 1}
\nu(U^i\cap W^j)\nu(V^k)
=1
\end{align*}

Thus,
$$ \liminf_{t\rightarrow \infty} \form{T^{\chi_U}_t (\chi_{V})}{\chi_{W}}
=\limsup _{t\rightarrow \infty} \form{T^{\chi_U}_t (\chi_V)}{\chi_W} =
\nu(U\cap W)\nu(V) $$

\end{proof}

\section{Proofs}

We are now in a position to gather the various ingredients supplied
in the last two sections into proofs.
Schematically, the sequence of proofs goes as follows:
$$ \mbox{prop~\ref{p:bounded} } + \mbox{ prop~\ref{p:limits} }
\Rightarrow
\mbox{ thm~\ref{operators} }
\Rightarrow \mbox{ thm~\ref{main} }
\Rightarrow \mbox{ thm~\ref{irr} } + \mbox{ thm~\ref{rigidity} }
$$
For a clarification of some of the notation, we advise the reader to
read Appendix~\ref{app-op}.

\begin{proof}[{Proof of theorem~\ref{operators}}]
The proof breaks into two parts.
In the first part, using the results of section~\ref{s:bounded}, we
prove that the maps $T_t$ have limit points when $t$ tends to infinity.
In the second part, using the results of section~\ref{s:limits},
we compute the actual limit.

For the first part, observe that the space of $\End(H)$
valued measures on $B$ can be naturally identified with
the space $(C(B)\otimes_{\pi} H \otimes_{\pi} H^*)^*$
(see corollary~\ref{op-valued-measure}).
By the Banach-Alaoglu theorem, every ball is compact in the corresponding
weak$^*$-topology.
Thus, we only have to prove that there is a uniform bound on the norms of
$T_t$.
The norms of $T_t$ are the operator norms of $T_t^1$.
Observe that $T_t^1$ preserve the subspace $L^{\infty}(B)<L^2(B)$.
Consider $T_t^1$ as transformations of $L^{\infty}(B)$,
using the positivity of all their coefficients as elements of $\BC\Gamma$,
their norms are given by $\|T_t^1(1)\|_{\infty}$,
which are uniformly bounded by proposition~\ref{p:bounded}.
By their self adjointness, the operators $T_t^1$ are uniformly bounded
as operators of $L^1(B)$ as well.
It follows, by the Riesz-Thorin interpolation theorem, that they
are uniformly bounded on $L^2(B)$,
and the first part of the proof follows.

We define two transformations, $p$ and $e$,
ranging from the algebraic tensor product $H\otimes H$, by
$$ p:H\otimes H \rightarrow \End(H)^*,\quad  p(u\otimes v)(T)=\form{Tu}{\bar{v}} $$
and
$$ e:H \otimes H \rightarrow \End(H), \quad e(u\otimes v)(w)=\form{w}{\bar{u}}v $$
These transformations and the relations between them are discussed in
Appendix~\ref{app-op}.

For the second part of the proof we choose a limit point at infinity
of $T_t$, denoted $T_{\infty}$,
and we wish to show that for every measurable set $U \subset B$,
$T_{\infty}(\chi_U)=e(1 \otimes \chi_U)$.
It is enough to prove this equation for sets $U$ in a generating set of the Borel
$\sigma$-algebra.
We will do so for the sets with measure zero boundaries.
Indeed, for such a set $U \subset B$, and for every two sets $W,V \subset B$
with measure zero boundaries,
by proposition~\ref{p:limits},

\begin{align*}
\form{p(\chi_V\otimes \chi_W)}{T_{\infty}^{\chi_U}}=
\form{T_{\infty}^{\chi_U} (\chi_V)}{\chi_W} = \nu(U\cap W)\nu(V) =
\\
\form{e(1\otimes \chi_U)\chi_V}{\chi_W} =
\form{p(\chi_V\otimes\chi_W)}{e(1\otimes\chi_U)}.
\end{align*}
By part 3 of lemma~\ref{>0},
the span of the set $\{p(\chi_V\otimes \chi_W)~|~\nu(\partial V)=\nu(\partial W)=0\}$
is weak$^*$ dense in $\End(L^2(B))^*$,
and the proof is complete.

\end{proof}

\begin{proof}[{Proof of theorem~\ref{main}}]
By theorem~\ref{operators},
$$ \{e(1\otimes\chi_U)~|~U \mbox{ measurable in } B\}\subset \VN^+(\rho),$$
As the group $\rho(\Gamma)$ is stable under conjugation, so is $\VN^+(\rho)$.
Then for every measurable $U,V \subset B$,
$$ e(\chi_U \otimes \chi_V)=e(1\otimes \chi_V)\circ e(1\otimes \chi_U)^*\in \VN^+(\rho).$$
We are done by part 4 of lemma~\ref{>0}.
\end{proof}

The implication to theorem~\ref{irr} is immediate.

\begin{proof}[Proof of theorem~\ref{irr}]
It follows at once from theorem~\ref{main} that $\VN(\rho)=\End(H)$.
This implies that the centralizer of $\rho(\Gamma)$ is the center of $\End(H)$,
that is, the group of scalar multiplications.
The irreducibility follows by Schur's lemma.
\end{proof}

In order to prove theorem~\ref{rigidity}
we state the following theorem, which in its
generality is due to Furman \cite{Furman-rigidityII}.

\begin{thm} \label{Furman}
Assume given a second Riemannian, negatively curved manifold $(M',g')$, with $\pi_1(M')\simeq \Gamma$.
Denote by $B'$ the boundary of the universal cover of $M'$,
and endow it with the corresponding Hausdorff measure (associated to some point in $\tilde{M}$).
The following are equivalent. \\
(1) There is a measurable $\Gamma$-equivariant measure class preserving
isomorphism of the Lebesgue spaces
$B$ and $B'$. \\
(2) The manifolds $(M,g)$ and $(M',g')$ give rise to
proportional Marked Length Spectra.
\end{thm}

\begin{proof}[{Proof of theorem~\ref{rigidity}}]
Assume given a Riemannian, negatively curved manifold $(M',g')$, with $\pi_1(M')\simeq \Gamma$.
Denote by $B'$ the boundary of the universal cover of $M'$,
and endow it with the corresponding Hausdorff measure (associated to some point in $\tilde{M}$).
Denote $H'=L^2(B')$, and by
$\rho':\Gamma\rightarrow \mathcal{U}(H')$,
the associated unitary representation.
Assume there is a non-trivial equivariant unitary intertwiner map
$L:H \rightarrow H'$.
By theorem~\ref{irr}, $L$ is an unitary representations isomorphism.
By equivarincy, a conjugation by $L$, gives a bijection
$\ad(L):\VN(\rho)^+ \simeq \VN(\rho')^+$.
Combining with theorem~\ref{main} we get
$\ad(L):\End(H)^+ \simeq \End(H')$.
By lemma~\ref{spec},
there exists an equivariant Lebesgue spaces isomorphism
$\ad(L)_*:B \rightarrow B'$.
We are done by theorem~\ref{Furman}.

\end{proof}

\appendix

\section{Integration on metric measure spaces} \label{app-mms}

In this appendix we recall the definition and some elementary properties of
metric $q$-regular measure spaces.
Through out this appendix let $(B,\metric)$ be a proper locally compact and complete metric space.
The space of Radon (i.e - locally finite and regular) measures on $B$ is denoted $M(B)$.
It is naturally embedded in $C_c(B)^*$, and we consider it with the (pull back of the) weak$^*$ topology.
The ball in $B$ of center $b$ and radius $t$ is denoted $B(b,t)$.
The following lemma is a clear consequence of the (inner) regularity assumption.

\begin{lemma} \label{l:lower-cont}
The real valued function $M(B) \times B\times [0,\infty) \rightarrow \BR$, $(\nu,b,r) \mapsto \nu(B(b,r))$
is lower semi-continuous.
\end{lemma}

\begin{defn} \label{d:mms}
A triple $(B,\metric,\nu)$ where $\nu \in M(B)$ is called a {\em metric measure space}.
Let $\eta$ be a positive number.
$(B,\metric,\nu)$ is called {\em $\eta$-regular} if there exist constants $0<k<k'$
such that for every $b\in B$ and $0<t\leq \diam(B)$,
\begin{equation} \tag{reg}
kt^{\eta} \leq \nu(B(b,t)) \leq k't^{\eta}.
\end{equation}
\end{defn}

\begin{lemma}
Assume $(B,\metric,\nu)$ is $\eta$-regular with multiplicative constants $0<k\leq k'$.
Fix numbers $0<s<t\leq \diam(B)$.
Let $f$ be a positive decreasing function on the interval $[s^{\eta},t^{\eta}]$. \\
Then for every $b\in B$,
\begin{align*}
& k \int_{s^{\eta}}^{t^{\eta}}f(u)du - (k'-k) f(s^{\eta}) \leq
  \int_{B(b,t)-B(b,s)} f(\metric(b,c)^{\eta}) d\nu(c) \\
& \leq k' \int_{s^{\eta}}^{t^{\eta}}f(u)du + (k'-k)s^{\eta} f(s^{\eta})
\end{align*}
\end{lemma}

\begin{proof}
We prove the right inequality, the proof of the left one being similar.
Fix $\epsilon>0$.
Choose a partition $u_0=s<u_1<\ldots<u_n=t$ of the interval $[s,t]$,
such that the upper sum associated to the partition $u_0^{\eta}<u_1^{\eta}<\ldots<u_n^{\eta}$
is less than $\int_{s^{\eta}}^{t^{\eta}}f(u)du+\epsilon$.
We simplify the notation by setting $v_i=f(u_{i-1}^{\eta})$, $|c|=\metric(b,c)$ and $B(u)=B(b,u)$.
We have
\begin{align*}
&   \int_{B(t)-B(s)} f(|c|^{\eta}) d\nu(c) =
   \sum_{i=1}^n \int_{B(u_i)-B(u_{i-1})} f(|c|^{\eta}) d\nu(c) \leq \\
&   \sum_{i=1}^n \int_{B(u_i)-B(u_{i-1})} v_i d\nu(c) =
   -v_1 \nu(B(u_0)) + \sum_{i=1}^{n-1} (v_i-v_{i+1}) \nu(B(u_i)) \\
&   + v_n \nu(B(u_n)) \leq
   -v_1 ku_0^{\eta} + \sum_{i=1}^{n-1} (v_i-v_{i+1}) k'u_i^{\eta} + v_n k'u_n^k = \\
&   v_1 (k'-k) u_0^{\eta} + k' \sum_{i=1}^n v_i (u_i^{\eta}-u_{i-1}^{\eta})  <
    (k'-k)s^{\eta} f(s) + k' \int_{s^{\eta}}^{t^{\eta}}f(u)du + \epsilon
\end{align*}
As $\epsilon$ was arbitrary we get the desired inequality.
\end{proof}

A particular case of interest is recorded in the following corollary.

\begin{cor} \label{int-as-log}
Assume that $(B,\metric,\nu)$ is $\eta$-regular with multiplicative constants $0<k\leq k'$,
and of diameter one.
Then for every $b\in B$ and $0<s<1$ we have
$$ -k\log s -(k'-k) \leq \int_{B(b,s)^c} \metric(b,c)^{-\eta} d\nu(c) \leq -k' \log s + (k'-k)$$
\end{cor}

In a regular space one can estimate an integral by sampling.
In order to make this statement precise we define what we mean by a "sample".

\begin{defn} \label{d:sampling}
Fix $r>0$.
An $r$-sampling set for $(B,\metric)$, with multiplicity $m\in \BN$,
is a finite set $S$ together with a map $\tilde{\ell}:S\rightarrow B$, $s\mapsto \tilde{\ell}_s$,
such that
$B=\cup_S B(\tilde{\ell}_s, r)$,
and for all $b\in B$, $|\tilde{\ell}^{-1}(B(b, r))| \leq m$.
\end{defn}

\begin{lemma} \label{l:sampling}
Assume that $(B,\metric,\nu)$ is $\eta$-regular with multiplicative constants $0<k\leq k'$.
Fix $r>0$.
Let $S$ be an $r$-sampling set for $(B,\metric)$, with multiplicity $m\in \BN$.
For every $L>0$, set $C_L=\frac{m (Le^L+1)k'}{k}$
(this constant does not depend on $r$).
Then for every function $f$ on $B$ such that $\log f$ is {\it $(r,L)$-almost continuous},
i.e., $|\log(f(x))-\log(f(y))|\leq L$ whenever $\metric(x,y)\leq r$, we have
$$  C_L^{-1} \frac{1}{|S|}\sum_S f(\tilde{\ell}_s)
\leq \int_B f(b)d\nu(b)
\leq C_L \frac{1}{|S|}\sum_S f(\tilde{\ell}_s). $$
\end{lemma}

\begin{proof}
We prove the right inequality, the proof of the left one being similar.
First observe that, denoting $g=\log f$, for every $b\in B$ and $b'\in B(b, r)$,
by almost continuity we have that
$$ |g(b)-g(b')| \leq  L, $$
hence also
$$ \max\{g(b),g(b')\} \leq g(b')+|g(b)-g(b')| \leq g(b')+  L. $$
Using the fact that for every $s,t\in\R$,
$$ |e^s-e^t|\leq |s-t|e^{\max\{s,t\}}, $$
we get
\begin{align*}
  f(b) = e^{g(b)} &
  \leq e^{g(b')} + |g(b)-g(b')|e^{\max\{g(b),g(b')\}} \\
& \leq e^{g(b')} + Le^{g(b')+L} = (Le^L+1)f(b').
\end{align*}
Second, observe that
$$ |S|\cdot k'r^{\eta}|\nu| =  \sum_S k'r^{\eta} \geq \int_B \sum_S \chi_{B(\ell_s, r)}(b) d\nu(b) \geq
\int_B d\nu = |\nu|, $$
therefore, $|S| r^{\eta} \geq  \frac{1}{k'}$.

Now, we are able to conclude that
\begin{align*}
& \frac{1}{|S|}\sum_S f(\tilde{\ell}_s) =
  \frac{1}{|S|}\sum_S \frac{1}{\nu(B(\tilde{\ell}_s, r))}\int_{B(\tilde{\ell}_s, r)} f(\tilde{\ell}_s) d\nu(b)\leq \\
& \frac{1}{|S|\nu(B(\tilde{\ell}_s, r))} \int_B \sum_S \chi_{B(\tilde{\ell}_s, r)} (Le^L+1)f(b) d\nu(b) \leq \\
& \frac{Le^L+1}{|S|kr^{\eta}} \int_B m  f(b) d\nu(b) \leq
  \frac{m (Le^L+1)k'}{k} \int_B  f(b) d\nu(b)
\end{align*}

\end{proof}

\section{Some functional analysis} \label{app-op}

Let $U$ and $V$ be Banach spaces (over $\BC$).
We denote by $\Hom(U,V)$ the space of continuous linear maps from $U$ to $V$.
In particular, we denote $U^*=\Hom(U,\BC)$.
We denote by $U\otimes V$ the algebraic tensor product of $U$ and $V$ over $\BC$.
Two interesting maps are given by
$$ p:U\otimes V \rightarrow \Hom(U,V^*)^*,\quad  p(u\otimes v)(T)=\form{Tu}{v} $$
and
$$ e:U \otimes V \rightarrow \Hom(U^*,V), \quad e(u\otimes v)(u^*)=\form{u^*}{u}v $$
The closures of the images are denoted, correspondingly, by
$U\otimes_{\pi} V$ and $U \otimes_{\varepsilon} V$.

Denote the composition of the maps
$$ \Hom(U,V^*)\overset{i}{\hookrightarrow} \Hom(U,V^*)^{**} \overset{p^*}{\rightarrow}
(U\otimes_{\pi} V)^* $$
by $p_*=p^*\circ i$.

\begin{lemma} \label{p_*}
The map $p_*$ is an isomorphism isometry of Banach spaces.
\end{lemma}

This is a well known fact. For completeness we give a brief of the proof.

\begin{proof}
Given an element $\lambda\in (U\otimes_{\pi}V)^*$ we consider the associated functional on the
{\em algebraic dual} of $U\otimes V$.
We, thus, obtain an (algebraic) map $\tilde{\lambda}$ from $U$ to the (algebraic) dual of $V$, given by
$\form{\tilde{\lambda}(u)}{v}=\lambda(u\otimes v)$.
It is readily checked that in fact $\tilde{\lambda}(U)\subset V^{*}$,
that the resulting operator $\tilde{\lambda}:U\rightarrow V^*$ is bounded by $\|\lambda\|_{\pi}$,
and that the map $\lambda\mapsto \tilde{\lambda}$ form an inverse map to $p_*$.
Thus $p_*$ is indeed a bijection, of norm one, and its inverse is of norm one.
\end{proof}

Let $(B,\metric,\nu)$ be a compact metric measure space.
Let $H=L^2(B,\nu)$ (over $\BC$).
Applying the lemma twice, we get

\begin{cor} \label{op-valued-measure}
The space of $\End(H)$-valued measure can be identified
with
$C(B)\otimes_{\pi} H \otimes_{\pi} H^*$.
\end{cor}

Set
$$ H \supset H^+=\{ u~|~u \geq 0 \mbox{ a.e}\}, \quad
H^* \supset (H^*)^+=\{ \lambda \in H^*~|~\lambda(H^+) \subset \BR_+ \}$$
$$ \mbox{and } \quad \End(H) \supset \End(H)^+=\{ T~|~TH^+ \subset H^+ \} $$
Recall that a subset of a topological vector space is called
{\em total} if its span is dense.

\begin{lemma} \label{>0}
Denote $E=\{ \chi_U~|~\nu(\partial U)=0\}$, and
$E\otimes E=
\{\chi_U \otimes \chi_V~|~\nu(\partial U)=\nu(\partial V)=0\}$.
\\
$1.$ The set $E$ is total in $H$. \\
$2.$
The closed cone in $H$ generated by
$E$ is $H^+$. \\
$3.$
The set $p(E\otimes E)$
is weak$^*$ total in $\End(H)^*$. \\
$4.$
The weak$^*$-closed cone in $\End(H)$, generated by
$e(E\otimes E)$
equals $\End(H)^+$.
\end{lemma}

\begin{proof}
$1.$ This follows from the regularity of $\nu$.
Indeed, for every measurable subset $A$ in $B$,
and for every $\epsilon$,
we can find compact set $K$ and open set $O$,
with $K \subset A \subset O$,
and $\nu(O-A),\nu(A-K) < \epsilon/2$.
All, but at most countably many of the sets
$$ K(t)=\{b \in B ~|~\metric(b,K) <t\}, \quad 0<t<\metric(K,O^c) $$
have zero measure boundary, and all of them
have symmetric difference with $A$ of measure less then $\epsilon$,
Thus,
$$ \overline{\spn\{ \chi_U~|~\nu(\partial U)=0\}}
= \overline{\spn\{ \chi_U~|~U \mbox{ measurable}\}}
= H
$$
$2.$ Every function in $h\in H^+$ can be approximated by elements of
$\spn(E)$.
But, every element of $\spn(E)$ is measurable with respect to some finite
($\sigma$-)algebra, the space of such  makes a (finite dimensional)
subspace $I< \spn(E)$,
and
the best approximation of $h$ in $I$ is given by its so called
conditional expectation, which is obviously positive. \\
$3.$
This follows from $1$,
using the fact that for every $u,v\in H$,
$\|u \otimes_{\pi} v\|=\|u\|\|v\|$. \\
$4.$
Fix $P \in \End(H)^+$.
By 1, a basis for the neighborhoods of $P$ is given by the sets
$$ U_{h_1,\ldots,h_n;\epsilon}=
\{ L~|~\mbox{for all } i=1,\ldots,n,~\|Ph_i-Lh_i\|<\epsilon \} \subset \End(H) $$
where  $\epsilon>0$, and the functions $h_1,\ldots,h_n$
are mutually orthogonal characteristic functions in $E$.
Assume given such a basis set.
Let $I= \spn\{h_1,\ldots,h_n\} < \spn(E)$.
Using 2, we can approximate the functions $Ph_1,\ldots,Ph_n$
by positive elements of $\spn(E)$, say $f_1,\ldots,f_n$.
There is a minimal finite ($\sigma$-)algebra, which measures all the $f_i$'s.
Denote by $h'_1,\ldots,h'_m \in E$ the characteristic functions of the atoms of that
algebra.
Set $I' =\spn\{h'_1,\ldots,h'_m\} < \spn(E)$.
Denote by $i:I \subset H$ and $i':I' \subset H$ the inclusions maps.
Let $L\in \Hom(I,I')$ be the transformation satisfies $Lh_i=f_i$.
Then $i\circ L\circ (i')^*$ is easily seen to be an element
of the positive cone generated by $e(E \otimes E)$,
which is in
$U_{h_1,\ldots,h_n;\epsilon}$
(for well enough approximating $f_i$'s).
\end{proof}

Finally, we have

\begin{lemma} \label{spec}
The iso-functor from Lebesgue spaces,
with measure class preserving maps,
to ordered von-Neumann algebras,
$$(B,\nu) \mapsto (\End(L^2(B,\nu)),\End(L^2(B,\nu))^+)$$
forms an equivalence of categories, with its image.
\end{lemma}

\begin{proof}
We will show how to reconstruct the Lebesgue space $B$,
given the couple $(\End(L^2(B,\nu)),\End(L^2(B,\nu))^+)$.
Observe that the set
$$\{ p \in \End(L^2(B,\nu)) ~|~ \mbox{both } p \mbox{ and }
1-p \mbox{ are orthogonal, positive projections}\}$$
coincides with the set
$$ \{ \chi_U ~|~U \mbox{ measurable in } B\} \subset L^{\infty}(B,\nu) < \End(L^2(B,\nu)) $$
Thus, the weak$^*$ closure of its span equals the commutative sub-von-Neumann
algebra $L^{\infty}(B,\nu)$,
which spectrum is isomorphic to $(B,\nu)$ as a Lebesgue space.
\end{proof}

\section{A "geodesics counting" formula} \label{app-Margulis}

We wish to express our deep gratitude to Fran\c{c}ois Ledrappier,
for pointing out to us the connection between the expression
$$ \frac{1}{|S_t|}
\sum_{\ga\in S_t} \chi_{\{\ga~|~\tilde{\ell}_{\ga^{-1}}\in U\}}
\chi_{\{\ga~|~\tilde{\ell}_{\ga} \in V\}}, $$
and the number of closed geodesics of a certain kind,
as well as the exact reference to its asymptotic,
computed in Margulis' thesis,
and by that enable us to replace a long argument by a short and elegant one.
The price we are willing to pay,
is having our presentation considerably less self contained.
This section is devoted to the explanation of
Ledrappier's observation.

The $\eta$-conformal density $\nu \rightarrow \Meas(B)$ can be naturally
interpreted as a collection of measures on the fibers of the map $T^1M \rightarrow M$,
by associating to every $m \in M$ the
measure $\nu_m$ on $T^1mM$ given by the
push forward of the measure $\nu_{\tilde{m}}$ under the
map
$$ B \rightarrow T^1_mM, \quad b \mapsto
dpr_{\tilde{m}}\left( \left. \frac{d\ell_{\tilde{m},b}(t)}{dt} \right|_{t=1} \right)
$$
Where $\tilde{m}$ is any pre-image of $m$ under the projection map $pr:X \rightarrow M$
(this definition does not depend on the choice of $\tilde{m}$).

Given two points $m,m'\in M$, two sets
$U\in T^1_mM$ and $U' \in T^1_{m'}M$ and an interval $I \subset \BR_+$,
we denote by $n(U,U',I)$ the number of geodesic segments $\ell:[0,s]\rightarrow M$ with
\begin{itemize}
\item $s\in I$.
\item $\ell(0)=m$ and $ \left. \frac{d\ell(t)}{dt} \right|_{t=0} \in -U$.
\item $\ell(s)=m'$ and $ \left. \frac{d\ell(t)}{dt} \right|_{t=s} \in U'$.
\end{itemize}


We refer the reader to theorem 7 in Margulis' thesis \cite[p. 57, Theorem 7]{Margulis-thesis}.
In a similar fashion to the proof of Theorem (6.2.5) in \cite{Yue},
applying Margulis' theorem to the submanifolds $T^1_mM$ and $T^1_{m'}M$
(which are transversal to both the stable and unstable foliations),
yields the following

\begin{thm}[Margulis]
Let $m,m' \in M$ be points,
$U \subset T^1_m$ and $U' \subset T^1_{m'}M$ measurable subsets satisfying
$\nu_m(\partial U)=\nu_{m'}(\partial U')=0$.
For every $a>0$, there exists $C=C_{a,m,m'}\in (0,\infty)$, such that
$$ \lim_{t\rightarrow \infty} e^{-\eta t}n(U,U',(t-a,t+a))=C\nu_m(U)\nu_{m'}(U') $$
(In particular, $C$ does not depend on the sets $U$ and $U'$).
\end{thm}

Observe that for the special case $pr(p)=m=m'$,
$$ n(U,U',I) = \sum_{\substack{\ga\in\Gamma \\ |\ga| \in I}}
\chi_U(\tilde{\ell}_{\ga^{-1}})\chi_{U'}(\tilde{\ell}_{\ga})
$$

Taking into account the definition of $\nu_m$ by means of $\nu$,
the fact that $\nu$ is normalized,
and using $S_t=\{\ga\in\Gamma~|~t-R < |\ga| < t+R\}$, we get

\begin{cor} \label{c:margulis}
For every $U,U'$ measurable subsets of $B$ with
$\nu(\partial U)=\nu(\partial U')=0$,
$$ \lim_{t\rightarrow \infty}
\frac{1}{|S_t|}
\sum_{\ga\in S_t}
\chi_U(\tilde{\ell}_{\ga^{-1}})\chi_{U'}(\tilde{\ell}_{\ga})
= \nu(U)\nu(U') $$
\end{cor}

\providecommand{\bysame}{\leavevmode\hbox to3em{\hrulefill}\thinspace}
\providecommand{\MR}{\relax\ifhmode\unskip\space\fi MR }
\providecommand{\MRhref}[2]{%
  \href{http://www.ams.org/mathscinet-getitem?mr=#1}{#2}
}
\providecommand{\href}[2]{#2}


\end{document}